\documentclass[12pt]{article}

\usepackage[T2A]{fontenc}
\usepackage[utf8]{inputenc}

\usepackage{amsfonts}
\usepackage{amssymb}
\usepackage{amsmath}
\usepackage{amsthm}

\def \le {\leqslant}
\def \ge {\geqslant}

\theoremstyle{plain}

\begin{document}
\begin{Huge}
\centerline{\bf Contribution to Vojt\v{e}ch Jarn\'{\i}k}
\end{Huge}
\begin{Large}
\vskip+0.5cm
\centerline{by Nikolay  Moshchevitin\footnote{ Research is supported by RFBR grant № 09-01-00371a and by Scientific School grant
   НШ-691.2008.1}}
\end{Large}
\vskip+2.0cm
\centerline{\bf Abstract}

We    find new inequalities between uniform and individual Diophantine exponents for three-dimensional Diophantine approximations.
Also we give a result for two linear forms in two variables. The results improves V.Jarn\'{\i}k's theorem (1954).
\vskip+4.0cm

{\bf 1. V.Jarn\'{\i}k's theorem.}
\,\,
This paper deals with real numbers only.
 We consider a matrix
$$
\Theta =\left(
\begin{array}{ccc}
\theta_{1}^1&\cdots&\theta_{1}^m\cr
\cdots&\cdots&\cdots\cr
\theta_{n}^1&\cdots&\theta_{n}^m
\end{array}
\right),
$$
and a system of linear forms
$$
 L_j ({\bf x})= \sum_{i=1}^m \theta_{j}^i x_i,\,\,\,
{\bf x}=(x_1,...,x_m).
$$
Let
$$
\psi_\Theta ( t)  = \,\,\,\,\,\min_{{\bf x}\in\mathbb{Z}^m: \, 0<M({\bf x})\le t
}\,\,\,\,\,\max_{1\le j\le n}||L_j({\bf x})||
 ,\,\,\, M({\bf x}) =\max_{1\le j\le m}|x_i|
$$
Suppose that $\psi_\Theta ( t) >0$ for all $t>0$ and define the
{\it uniform} Diophantine exponent $\alpha (\Theta)$ as the supremum of the set
$$\{\gamma >0:\,\,\,
\limsup_{t\to +\infty} t^\gamma
\psi_\Theta (t) <+\infty \},
$$
 The
{\it individual} Diophantine exponent $\beta (\Theta)$ is defined as the supremum of the set
$$
\{\gamma >0:\,\,\,
\liminf_{t\to +\infty} t^\gamma
\psi_\Theta (t) <+\infty\}
.
$$
In \cite{JRUS} V.Jarn\'{\i}k proved the following result.

{\bf Theorem 1.}\,\,{\it
Suppose that $\psi_\Theta ( t) >0$ for all $t>0$.

{\rm (i)} Consider the case $m=1$ and suppose that among numbers $\theta^1_1,...,\theta^1_n$ there exist at least two numbers
which are linearly independent together with $1$ over $\mathbb{Z}$. Then
\begin{equation}
\label{1}
\beta(\Theta) \ge \frac{(\alpha(\Theta))^2}{1-\alpha (\Theta)}.
\end{equation}

{\rm (ii)} Consider the case $m=2$. Then
\begin{equation}
\label{2}
\beta(\Theta) \ge \alpha(\Theta)(\alpha(\Theta) -1).
\end{equation}

{\rm (iii)} Consider the case $m>2$.
Suppose that $\alpha (\Theta) >(5m^2)^{m-1}$.
 Then
\begin{equation}
\label{3}
\beta(\Theta) \ge (\alpha(\Theta))^{\frac{m}{m-1}} - 3\alpha(\Theta).
\end{equation}
}

It is a well-known fact that under the conditions of Theorem 1 one has
 $$
\frac{m}{n}\le\alpha (\Theta)\le \beta (\Theta)\le +\infty.
$$
Moreover for $m=1$ one has
$$
\frac{1}{n}\le \alpha (\Theta)\le 1.$$
In \cite{L}
M.Laurent proved that in the cases $m=1,n=2$ and $m=2,n=1$ the inequalities of V.Jarn\'{\i}k's theorem cannot be improved.

In the
present note we give a sketched proof of an improvement of V.Jarn\'{\i}k's theorem  in the cases
$m=1,n=3$, $m=3,n=1$ and $m=n=2$.

{\bf 2. Results.}
For $\alpha >0$ put
$$
g_1(\alpha ) =
\frac{\alpha (1-\alpha) +\sqrt{\alpha^2(1-\alpha)^2 +4\alpha (2\alpha^2-2\alpha+1)}}{2(2\alpha^2-2\alpha+1)}.
$$
The value $g_1(\alpha)$ is the largest root of the equation
$$
(2\alpha^2-2\alpha +1)x^2+\alpha(\alpha -1)x -\alpha = 0.$$
Given $\alpha $ consider a system of equations
\begin{equation}\label{sys}
\gamma = \frac{1}{\alpha}+\frac{\alpha -1}{\alpha} \frac{\delta}{\gamma}=
\frac{\alpha}{\gamma(1-\alpha)-\alpha}.
\end{equation}
Then there exist a solution of this system with $\delta=g_1(\alpha)$.
Note that
$$
g_1(1/3)=g_1(1)=1
$$
and for $1/3<\alpha <1$ one has $g_1(\alpha )>1$.
Let $\alpha_0$ be the unique real root of the equation
$$
x^3-x^2+2x - 1 =0.
$$
In the interval $ 1/3< \alpha <\alpha_0$
one has
$$
g_1(\alpha ) > \max \left(1, \frac{\alpha}{1-\alpha}\right).
$$

{\bf Theorem 2.}\,\,{\it Suppose that $m=1,n=3$ and the matrix $
\Theta=\left(\begin{array}{c}
\theta_1\cr \theta_2 \cr \theta_3\end{array}\right)$ consists of numbers
linearly independent over   $\mathbb{Z}$  together with $1$. Then
\begin{equation}\label{novoe0}
\beta (\Theta) \ge \alpha (\Theta) g_1 (\alpha (\Theta)).
\end{equation}
}

The inequality (\ref{novoe0}) is better than (\ref{1}) in the interval $1/3 < \alpha(\Theta )<\alpha_0$.

For  $\alpha\ge 3$ define
$$
g_2(\alpha ) =\sqrt{\alpha+\frac{1}{\alpha^2}-\frac{7}{4}}+\frac{1}{\alpha}-\frac{1}{2},\,\,\,\,
h(\alpha ) = \alpha - g_2(\alpha ) -1.
$$
Here we should note that the functions
   $g_2(\alpha) $  and  $h(\alpha)$
monotonically increase to $+\infty$ when
  $\alpha \to+\infty$   and
$$
g_2(3) = h(3) = 1,\,\,\,\,  g_2(\alpha ) \le \alpha - 2.
$$

{\bf Theorem 3.}\,\,{\it Suppose that $m=3,n=1$ and the matrix $\Theta = (\theta^1,\theta^2,\theta^3)$ consists of numbers
linearly independent over   $\mathbb{Z}$  together with $1$. Then
\begin{equation}\label{novoe}
\beta (\Theta) \ge \alpha (\Theta) g_2 (\alpha (\Theta)).
\end{equation}

 }

The inequality (\ref{novoe}) is better than (\ref{3}) for all values of $\alpha(\Theta)$.

For $\alpha \ge 1$
put
$$
g_3 (\alpha ) =
\frac{1-\alpha +\sqrt{(1-\alpha)^2+4\alpha (2\alpha^2-2\alpha+1)}}{2\alpha}
.
$$
So $g_3 (\alpha ) $
is a solution of the equation
\begin{equation}\label{eq}
\alpha x^2+(\alpha-1)x - (2\alpha^2-2\alpha+1) =0.
\end{equation}
We see that $g_3   (1) = 1$ and for $\alpha >1$ one has $g_3(\alpha ) >1$.
Moreover in the interval
$$1\le\alpha < \left(\frac{1+\sqrt{5}}{2}\right)^2
$$
one has
$$g_3 (\alpha ) >\max (1, \alpha - 1).
$$
{\bf Theorem 4. }\,\,{\it
 Consider four real numbers $\theta_j^i,\,\, i,j=1,2$ linearly independent over $\mathbb{Z}$ together with 1.
Let $m=n=2$ and consider the matrix
$$
\Theta =
\left(
\begin{array}{cc}
\theta^1_1&\theta^2_1\cr
\theta^1_2&\theta^2_2
\end{array}
\right)
$$
Then
\begin{equation}
\label{2mod}
\beta(\Theta) \ge \alpha(\Theta)g_3(\alpha(\Theta)).
\end{equation}
 }

Theorem 4 improves Theorem 1 for $\alpha (\Theta) \in \left(1,\left(\frac{1+\sqrt{5}}{2}\right)^2\right)$.

{\bf 3. Best approximations.}

For integer vector $ {\bf x} = (x_1,...,x_m) \in \mathbb{Z}^m$ put
$$
  \zeta ({\bf x}) =  \max_{1\le j\le n}||L_j({\bf x})||.
$$
A point
  $ {\bf x} = (x_1,...,x_m) $ is defined to be  a {\it best approximation} if
$$
\zeta ({\bf x})=\min_{{\bf x}'} \zeta ({\bf x}'),$$
where the minimum is taken over all  $ {\bf x}' = (x_1',...,x_m')
\in \mathbb{Z}^m $ such that
$$
0<  M({\bf x}')\le M({\bf x}). $$
(For each best approximation   $ {\bf x}$, the point
 $ -{\bf x}$ is also a best approximation.) Consider the case when all numbers
$\theta_{j}^i,\,\,\, 1\le i \le m,\,\,\, 1\le j\le n$  are linearly independent over  $\mathbb{Z}$ together with $1$.
Then all best approximations form the sequences
$$
 {\bf x}_1,{\bf x}_2,...,{\bf x}_\nu, {\bf x}_{\nu+1}, ...
$$
 $$
M({\bf x}_1)< M({\bf x}_2)<...<M({\bf x}_\nu)<M({\bf x}_{\nu+1})<... ,
 $$
$$
\zeta({\bf x}_1)> \zeta({\bf x}_2)>...>\zeta({\bf x}_\nu)>\zeta({\bf x}_{\nu+1})>... . $$
We use the notation
$$
M_\nu = M({\bf x}_\nu),\,\,\,\,\,\zeta_\nu = \zeta({\bf x}_\nu).
$$
Define $ y_{1,\nu},...,y_{n,\nu}\in \mathbb{Z}^n$ to be integers such that
$$
||L_j({\bf x}_\nu)||= |L_j({\bf x}_\nu)+y_{j,\nu}|.
$$
 We need the notation
$$
{{\bf z}_\nu} = (x_{1,\nu},...,x_{m,\nu}, y_{1,\nu},..., y_{n,\nu} ) \in \mathbb{Z}^{d},\,\, d = m+n.
$$
If
$$
\psi_{\Theta}(t) \le \psi (t)$$
with continuous and decreasing to zero function $\psi (t)$
one can easily see that
$$
\zeta_\nu \le \psi (M_{\nu+1}).$$

 Some useful fact about best approximations one can find in \cite{ME}.

{\bf 4. Sketched proof of Theorem 2.}\,\,

Suppose that
$$
\psi_{\Theta}(t) \le \psi (t)$$
with some  continuous function $\psi(t)$ decreasing to zero as $t\to+\infty$.
Moreover we suppose that the function $t\mapsto t\cdot \psi (t)$ increases to infinity as $t\to +\infty$.

Consider best approximation vectors ${\bf z}_\nu = (x_\nu,y_{1,\nu},y_{2,\nu},y_{3,\nu})$.
From the condition that numbers $1,\theta_1,\theta_2,\theta_3$ are linearly independent over $\mathbb{Z}$ we see that
   there exist infinitely many pairs of indices
   $\nu<k, \nu\to +\infty$   such that

$\bullet $  both triples
$$
{\bf z}_{\nu-1},{\bf z}_\nu,{\bf z}_{\nu+1};\,\,\,\,\,\,
{\bf z}_{k-1},{\bf z}_k,{\bf z}_{k+1}
$$ consist of
  linearly independent vectors;

$\bullet$
there exists a two-dimensional linear subspace   $\pi$  such that
$$
{\bf z}_l\in \pi,\,\,\, \nu\le l\le k;\,\,\,\,\, {\bf z}_{\nu-1} \not\in \pi,\,\,\, {\bf z}_{k+1} \not\in \pi;
$$

$\bullet$  the vectors
$$
{\bf z}_{\nu-1},{\bf z}_\nu,{\bf z}_{k},{\bf z}_{k+1}
$$
are linearly independent.

So
 $$
1\le |
{\rm det}
\left(
\begin{array}{cccc}
y_{1,\nu-1}&y_{2,\nu-1}&y_{3,\nu-1}& x_{\nu-1}\cr
y_{1,\nu}&y_{2,\nu}&y_{3,\nu}& x_{\nu}\cr
y_{1,k}&y_{2,k}&y_{3,k}& x_{k}\cr
y_{1,k+1}&y_{2,k+1}&y_{3,k+1}& x_{k+1}
\end{array}
\right)|
\le
$$
\begin{equation}\label{promejutok0}
\le
24 \zeta_{\nu-1}\zeta_\nu\zeta_{k} M_{k+1} \le
24 \psi (M_\nu) \psi (M_{\nu +1})\psi (M_{k+1})M_{k+1} .
 \end{equation}
 Consider three cases.

{\bf 1$^0$.} Given $\gamma >\frac{\alpha (\Theta)}{1-\alpha (\Theta)}$ there exist infinitely many pairs  $(\nu,k)$
such that
$$
M_{k+1}\le M_{\nu+1}^{\gamma}.
$$
From  (\ref{promejutok0}) we deduce that
 $$
 \frac{1}{24\psi(M_\nu)}\le  M_{\nu+1}^\gamma\cdot \psi (M_{\nu+1})\cdot \psi (M_{\nu+1}^\gamma )
$$
Suppose the function $t\mapsto t^\gamma \cdot \psi (t)\cdot \psi (t^\gamma )$ to be increasing and let
$\rho (t)$ be the inverse function.
 We see that
\begin{equation}\label{a}
\zeta_\nu\le \psi (M_{\nu+1})\le
\psi\left(
\rho \left(
\frac{1}{24\psi (M_\nu)}
\right)
\right).
 \end{equation}

{\bf 2$^0$.} Given $\delta \ge 1$.  There exist infinitely many pairs  $(\nu,k)$
such that
$$
M_{k+1}\ge  M_k^{\delta}.
$$
Then we immediately have
\begin{equation}\label{b}
\zeta_k \le \psi (M_{k+1}) \le \psi ( M_k^{\delta}).
\end{equation}

{\bf 3$^0$.}  There exist infinitely many pairs  $(\nu,k)$
such that
$$
 M_{\nu+1}^\gamma \le M_{k+1}\le M_k^\delta.
$$

{\bf Lemma 1.}\,\,{\it
Let $ {\bf z}_\nu , a\le \nu \le b$  lie in a two-dimensional linear subspace $\pi \subset \mathbb{R}^4$. Then
for all $\nu_1,\nu_2$
from the interval $a\le \nu_j\le b-1$ one has
$$
\zeta_{\nu_1}M_{\nu_1+1} \asymp_\Theta \zeta_{\nu_2}M_{\nu_2+1}.
$$}

{\bf  Sketched proof of Lemma 1.}\,\,
Consider two-dimensional  lattice $\Lambda = \pi \cap \mathbb{Z}^4$.
The parallelepiped
$$
\{
{\bf z}=(x,y_1,y_2,y_3):\,\,
|x|< M_{k+1},\,\, \max_{1\le j\le 3}|\theta_jx-y_j|<\zeta_k\}
$$
has no non-zero integer points inside for every $k$.
So for $a\le \nu\le b-1$ we have
$$
\zeta_{\nu}M_{\nu+1} \asymp_\Theta {\rm det}_2\Lambda.
$$
Lemma 1 follows.

We apply Lemma 1 and obtain the inequality
\begin{equation}\label{c}
 \zeta_{k-1}\ll_\Theta
\frac{M_{\nu+1}\psi (M_{\nu+1})}{M_k}\ll_\Theta
M_k^{\frac{\delta}{\gamma}-1}\psi (M_k^{\frac{\delta}{\gamma}})\le
M_{k-1}^{\frac{\delta}{\gamma}-1}\psi (M_{k-1}^{\frac{\delta}{\gamma}})
\end{equation}
(here we suppose the function $t\mapsto
t^{\frac{\delta}{\gamma}-1}\psi (t^{\frac{\delta}{\gamma}})
 $ to be decreasing).

Given $\varepsilon >0$ we put $\psi(t) =t^{-\alpha+\varepsilon}$ to deduce (\ref{novoe0}) from (\ref{a},\ref{b},\ref{c}).
Here we should take into account that $g_1(\alpha)$ satisfies the system (\ref{sys}).

{\bf 5. Sketched proof of Theorem 3.}\,\,
Suppose that
$$
\psi_{\Theta}(t) \le \psi (t)$$
with some  continuous function $\psi(t)$ decreasing to zero as $t\to+\infty$.

First of all consider the case when there exists a 3-dimensional linear subspace $\Pi$ such that all vectors
${\bf z}_\nu =(x_{1,\nu},x_{2,\nu},x_{3,\nu},y_\nu)$
belong to $\Pi$ for all $\nu$ large enough. Then obviously we can apply the statement (ii) of Theorem 1 and (\ref{2}) gives the bound which is better than
  (\ref{novoe}).

So we can suppose that  there exist infinitely many pairs of indices
   $\nu<k,\nu\to +\infty$   such that

$\bullet $  both triples
$$
{\bf z}_{\nu-1},{\bf z}_\nu,{\bf z}_{\nu+1};\,\,\,\,\,\,
{\bf z}_{k-1},{\bf z}_k,{\bf z}_{k+1}
$$ consist of
  linearly independent vectors;

$\bullet$
there exists a two-dimensional linear subspace   $\pi$  such that
$$
{\bf z}_l\in \pi,\,\,\, \nu\le l\le k;\,\,\,\,\, {\bf z}_{\nu-1} \not\in \pi,\,\,\, {\bf z}_{k+1} \not\in \pi;
$$

$\bullet$  the vectors
$$
{\bf z}_{\nu-1},{\bf z}_\nu,{\bf z}_{\nu+1},{\bf z}_{k+1}
$$
are linearly independent.

Now we see that
\begin{equation}\label{promejutok}
1\le|
{\rm det}
\left(
\begin{array}{cccc}
x_{1,\nu-1}&x_{2,\nu-1}&x_{3,\nu-1}& y_{\nu-1}\cr
x_{1,\nu}&x_{2,\nu}&x_{3,\nu}& y_{\nu}\cr
x_{1,\nu+1}&x_{2,\nu+1}&x_{3,\nu+1}& y_{\nu+1}\cr
x_{1,k+1}&x_{2,k+1}&x_{3,k+1}& y_{k+1}
\end{array}
\right)|
\le
24 \zeta_{\nu-1}M_\nu M_{\nu+1}
M_{k+1} \le
24 \psi (M_\nu) M_\nu M_{\nu+1}M_{k+1}.
\end{equation}
 Consider three cases.

{\bf 1$^0$.}  There exist infinitely many pairs  $(\nu,k)$
such that
$$
M_{k+1}\le M_\nu^{h(\alpha(\Theta))}.
$$
From  (\ref{promejutok}) we deduce that
$$
M_{\nu+1} \ge \frac{1}{24\psi(M_\nu)M_\nu^{1+h(\alpha(\Theta))}},
\,\,\,\,
\zeta_\nu \le \psi\left(\frac{1}{24\psi(M_\nu)M_\nu^{1+h(\alpha(\Theta))}}
\right).
$$
The inequality (\ref{novoe}) follows  from the last inequality immediately.

{\bf 2$^0$.}  There exist infinitely many pairs  $(\nu,k)$
such that
$$
M_{k+1}\ge  M_k^{g_2(\alpha(\Theta))},
$$
Then we immediately have
$$
\zeta_k \le \psi (M_{k+1}) \le \psi ( M_k^{g_2(\alpha(\Theta))}),
$$
and  (\ref{novoe}) follows.

{\bf 3$^0$.}  There exist infinitely many pairs  $(\nu,k)$
such that
$$
M_\nu^{h(\alpha(\Theta))}\le
M_{k+1}\le  M_k^{g_2(\alpha(\Theta))}.
$$

{\bf Lemma 2.}\,\,{\it
Let $ {\bf z}_\nu , a\le \nu \le b$  lie in a two-dimensional linear subspace $\pi \subset \mathbb{R}^4$. Then
for all $\nu_1,\nu_2$
from the interval $a\le \nu_j\le b-1$ one has
$$
\zeta_{\nu_1}M_{\nu_1+1} \asymp_\Theta \zeta_{\nu_2}M_{\nu_2+1}.
$$}

{\bf  Sketched proof of Lemma 2.}\,\,
 Consider the projection of the subspace  $\pi$
onto the $3$-dimensional subsbace
   $${\cal L}(\Theta)
=\{{\bf z}=(x_1,x_2,x_3,y):\,\,
\theta^1x_1+\theta^2x_2+\theta^3x_3 +y =0\}.
$$
 In the general situation this projection is a two-dimensional subspace $\pi^*$. The intersection
   $\ell = \pi\cap \pi^*$ form a one-dimensional subspace.
The distance from a point ${\bf z}\in \pi$ to the subspace ${\cal L}$
is proportional to the distance from ${\bf z} $ to $\ell$.
Let   $\delta$  be the coefficient of this proportionality
 The vectors
  ${\bf z}_l$
become the vectors of the best approximations (associated with the induced norm on $\pi$) from a lattice $\Lambda =\mathbb{Z}^4\cap\pi
$
to the one-dimensional subspace $\ell$.
By the Minkowski convex body theorem applied to the two-dimensional lattice $\Lambda$ we deduce that
  $$
\gamma_1(\Theta
)  \delta \,{\rm det}\,\Lambda \le
\zeta_\nu M_{\nu+1} \le \gamma_2(\Theta
)  \delta \,{\rm det}\,\Lambda ,\,\,\,\,\,\,
a\le \nu\le b-1
$$
with some positive constants $\gamma_i(\Theta),\, i =1,2$ depending on $\Theta$.
 Lemma 2  is proved.

In the case {\bf 3}$^0$ we apply Lemma 2  an obtain the inequality
$$
\zeta_\nu \ll_\Theta
\frac{\psi (M_k)M_k}{M_{\nu+1}}\ll_\Theta
\psi\left(M_\nu^{\frac{h(\alpha(\Theta))}{g_2(\alpha(\Theta))}}\right)
M_\nu^{\frac{h(\alpha(\Theta))}{g_2(\alpha(\Theta))}-1}
.
$$
As
$$
\alpha (g_2(\alpha))^2 +(\alpha -2)g_2(\alpha) -(\alpha -1)^2 =0,
$$
we deduce (\ref{novoe}) in each case.

 {\bf 6. Sketched proof of Theorem 4.}\,\,

Define
$$
R(\Theta) =  \min_\Pi\, {\rm dim} \,
\Pi \ge 2,
$$
where the minimum is taken over all linear subspaces $\Pi \subseteq \mathbb{R}^4$ such that
there exists $\nu_0$ such that for all $\nu\ge \nu_0$ one has
${\bf z}_\nu =(x_{1,\nu},x_{2,\nu},y_{1\nu},y_{2,\nu})\in \Pi$.

We consider several cases.

{\bf 1$^0$.}
Suppose that $R(\Theta) = 2$. Then  (by Theorem 8 from \cite{ME}, Section 2.1, applied to the case $m=n=2$)
we see that $\alpha (\Theta ) = 1$ and there is nothing to prove.

{\bf 2$^0$.} The equality $R(\Theta) =3$ is not possible (see Corollary 4 from the Section 2.1 from \cite{ME}).

{\bf 3$^0$.}
Suppose that $R(\Theta) = 4$.
Then
   there exist infinitely many pairs of indices
   $\nu<k, \nu\to +\infty$   such that

$\bullet $  both triples
$$
{\bf z}_{\nu-1},{\bf z}_\nu,{\bf z}_{\nu+1};\,\,\,\,\,\,
{\bf z}_{k-1},{\bf z}_k,{\bf z}_{k+1}
$$ consist of
  linearly independent vectors;

$\bullet$
there exists a two-dimensional linear subspace   $\pi$  such that
$$
{\bf z}_l\in \pi,\,\,\, \nu\le l\le k;\,\,\,\,\, {\bf z}_{\nu-1} \not\in \pi,\,\,\, {\bf z}_{k+1} \not\in \pi;
$$

$\bullet$  the vectors
$$
{\bf z}_{\nu-1},{\bf z}_\nu,{\bf z}_{k},{\bf z}_{k+1}
$$
are linearly independent.

So
 $$
1\le |
{\rm det}
\left(
\begin{array}{cccc}
x_{1,\nu-1}&x_{2,\nu-1}&y_{1,\nu-1}& y_{2,\nu-1}\cr
x_{1,\nu}&x_{2,\nu}&y_{1,\nu}& y_{2,\nu}\cr
x_{1,k}&x_{2,k}&y_{1,k}& x_{2,k}\cr
x_{1,k+1}&x_{2,k+1}&y_{1,k+1}& y_{2,k+1}
\end{array}
\right)|
\le
24 \zeta_{\nu-1}\zeta_\nu M_k  M_{k+1}
.
$$
Suppose that a function $\psi (t) $ decrease to zero as $t\to+\infty$. Suppose also that the function
$t\mapsto t\cdot \psi (t)$ also decrease to zero.
Suppose that
$$
\psi_\Theta (t) \le \psi (t)$$
for all positive $t$. Then
$\zeta_l \le \psi (M_{l+1}),\, l = 1,2,3,..$ and
\begin{equation}\label{uro}
1\le 24 M_{k+1}M_k\psi (M_{\nu+1})\psi (M_\nu).
\end{equation}
We must consider two subcases.

{\bf 3.1$^0$.}
 Given $\gamma >1$ there exist infinitely many pairs  $(\nu,k)$
such that
$$
M_{k+1}\ge  M_k^{\gamma}.
$$
Then we immediately have
\begin{equation}\label{b}
\zeta_k \le \psi (M_{k+1}) \le \psi ( M_k^{\gamma}).
\end{equation}

{\bf 3.2$^0$.}  There exist infinitely many pairs  $(\nu,k)$
such that
$$
  M_{k+1}\le M_k^\gamma.
$$
Then from (\ref{uro}) we see that
\begin{equation}\label{urod}
M_k \ge (\psi (M_\nu ))^{-\frac{2}{1+\gamma}}.
\end{equation}
 Consider two-dimensional  lattice $\Lambda = \pi \cap \mathbb{Z}^4$ with the fundamental two-dimensional volume
${\rm det }\,\Lambda$.

We may suppose that ${\rm dim}\, {\cal L}\cap \pi=1$
(the case ${\cal L}\cap \pi ={\bf 0}$
can be considered in a similar way).

For any point ${\bf z}\in \pi$ the distance from ${\bf z}$ to the two-dimensional  linear subspace
$$
{\cal L} =\{ {\bf z}=(x_1,x_2,y_1,y_2):\,\,
\theta_1^1x_1+\theta^2_1x_2+y_1 =
\theta_2^1x_1+\theta^2_2x_2+y_2=0\}
$$
is proportional to the distance from
${\bf z}$ to the one-dimensional  linear subspace
${\cal L}\cap \pi$.
Let $\delta$ be the coefficient of this proportionality (the "angle"   between two-dimensional subspaces
$\pi$    and ${\cal L}$).
  The parallelepiped
$$
\{
{\bf z}=(x_1,x_2,y_1,y_2):\,\,
|x|< M_{l+1},\,\, \max_{1\le j\le 3}|\theta_jx-y_j|<\zeta_l\}
$$
has no non-zero integer points inside for every $l$.
 From the Minkowski convex body theorem e wee that
  \begin{equation}\label{uroda}
\gamma_1(\Theta
)  \delta \,{\rm det}\,\Lambda \le
\zeta_l M_{l+1} \le \gamma_2(\Theta
)  \delta \,{\rm det}\,\Lambda ,\,\,\,\,\,\,
\nu\le l\le k-1
\end{equation}
with some positive constants $\gamma_i(\Theta),\, i =1,2$ dependind on $\Theta$.
 So from (\ref{urod}) and (\ref{uroda})  we see that
\begin{equation}\label{c}
 \zeta_\nu \ll_\Theta
\frac{\psi (M_k)M_k}{M_{\nu+1}}\ll_\Theta
M_\nu^{-1}(\psi (M_\nu))^{-\frac{2}{1+\gamma}}\psi\left( (\psi(M_\nu))^{-\frac{2}{1+\gamma}}\right).
\end{equation}
We should consider $ \psi (t) = t^{-\alpha(\Theta)+\varepsilon}$ for small positive $\varepsilon$. As
$\gamma =g_3(\alpha (\Theta)) $ satisfies (\ref{eq}) Theorem 2 follows.


\begin{thebibliography}{100}


\bibitem{JRUS}
В.Ярник,\,\,\,
К теории однородных линейных диофантовых приближений.
//
Чехословацкий математический журнал, т. 4 (79), 330 - 353 (1954).


\bibitem{L} M.Laurent,\,\,\,
Exponents of Diophantine approximations in dimension two.//
Canad.J.Math. 61, 1 (2009),165 - 189; preprint available at arXiv:math/0611352v1 (2006).


\bibitem{ME}
N.G.Moshchevitin,\,\,\,
Khintchine's singular systems and their applications.//Russian Math. Surveys.   65:3 43 - 126 (2010);(2010) Preprint available at arXiv:    (2009).

\end{thebibliography}
\end{document}